\documentclass[12pt]{article}
\usepackage{amssymb,amsmath,graphicx,amsthm,mathrsfs}
\newtheorem{Theorem}{Theorem}[section]
\newtheorem{Lemma}[Theorem]{Lemma}

\newtheorem{Corollary}[Theorem]{Corollary}

\newtheorem{Remark}[Theorem]{Remark}
\numberwithin{equation}{section}

 \def\to{\rightarrow}

\newcommand{\q}{\quad}    \def\R{{\mathbb R}}
  \def\f{\frac}

 \def\t{\times}

 \def\v{{\bf v}}

\newcommand{\lc}
{\mathrel{\raise2pt\hbox{${\mathop<\limits_{\raise1pt\hbox
{\mbox{$\sim$}}}}$}}}

\newcommand{\gc}
{\mathrel{\raise2pt\hbox{${\mathop>\limits_{\raise1pt\hbox{\mbox{$\sim$}}}}$}}}

\newcommand{\ec}
{\mathrel{\raise2pt\hbox{${\mathop=\limits_{\raise1pt\hbox{\mbox{$\sim$}}}}$}}}

\def\bb{\begin{equation}} \def\ee{\end{equation}}

\def\beqn{\begin{eqnarray}}  \def\eqn{\end{eqnarray}}

\def\bbx{\begin{equation*}}   \def\eex{\end{equation*}}

\def\beqnx{\begin{eqnarray*}} \def\eqnx{\end{eqnarray*}}

\def\bn{\begin{enumerate}} \def\en{\end{enumerate}}

\def\bd{\begin{description}} \def\ed{\end{description}}

\def \d  {\,\mathrm{d}}
\def \f  {\frac}

\def \q  {\quad}
\def \v {\varphi}
\def \lg {\langle}
\def \rg {\rangle}

\begin{document}
\title{The time optimal control
with  constraints of the rectangular type
for linear time-varying ODEs}
\author{Can Zhang\thanks{
School of Mathematics and Statistics, Wuhan University, Wuhan,
430072, China. (zhangcansx@163.com). The author was partially supported by the National Natural Science Foundation of China under grants
11161130003 and 11171264.} }

\date{}
\maketitle
\begin{abstract}
In this paper, we study a time optimal control
problem of some linear time-varying
ordinary differential equations, where the control constrained set is of the rectangular
type. We aim to build up a necessary and sufficient condition and provide  an algorithm
for the optimal time, as well as the optimal control.  We first set up a norm optimal control problem associated with the control constraints of the rectangular
type; then establish an equivalence theorem between the time  optimal control problem  and the aforementioned norm optimal control problem; finally, reach the aim, through utilizing
the equivalence theorem and analyzing the variational characterization for the norm optimal control problem.
\\

\noindent\textbf{Key words.} time optimal control, norm optimal control, optimal time, optimal norm, control constraints of the rectangular type\\

\noindent\textbf{AMS Subject Classifications.} 49J15, 49K15
 \end{abstract}

\section{Introduction}
$\;\;\;\;$ Let $m$ and $d$ be two natural numbers. Let $A(\cdot)\in C([0,+\infty);\R^{m\times m})$ and  $b_i\in \R^{m}$ with $i=1,\cdots,d$.  Consider the following  controlled
 linear  time-varying ordinary differential equation:
\begin{eqnarray}\label{1}
 y'(t)+A(t)y(t)=\sum_{i=1}^{d}b_iu^i(t),\q t\geq0,\q
 y(0)=y_0.
\end{eqnarray}
Here and throughout this paper, the initial state $y_0$ is assumed to be a nonzero vector in $\mathbb{R}^m$,  $u^i(\cdot)$, $i=1,\cdots,d$, are  control functions from $\mathbb{R}^+$ to $\mathbb{R}^1$. The following notations will be frequently used in this paper: we denote by  $y(\cdot;u)$ the solution of Equation (\ref{1}) corresponding to the control $u=(u^1, \cdots, u^d)$;
 write $\lg\cdot,\cdot\rg$ and $\|\cdot\|$  for the usual inner product
 and Euclid norm in $\R^m$ respectively; use
$A^*$ and $\|A\|_{\R^{m\times m}}$ to denote the transpose
 matrix and matrix norm of $A$ respectively.
 The following assumptions on  $A(\cdot)$ will be effective throughout the paper:
{
\it

\noindent (H.1) $A(\cdot)$ is real analytic on $(0,+\infty)$.

\noindent (H.2) For each $i=1,\cdots,d$, $(A(\cdot),b_i)$
satisfies  Conti's condition, namely, the equality:
$$
\displaystyle \int_0^{+\infty}|\lg b_i,\v(t)\rg|\d t=+\infty
$$
 holds for each nonzero solution to the  dual equation:
\bb \label{3}
\varphi'(t)-A^*(t)\varphi(t)=0,\;\;t\in [0,+\infty).
\ee
}
With respect to  Conti's condition, we refer the readers to  \cite{CON} or \cite{SB}.
It is worth to  mention that  Conti's condition holds for $A(\cdot)$ and $b_1\in \mathbb{R}^m$ if and only if the system (\ref{1}) where $d=1$ is null controllable with the control constraints $|u^1(t)|\leq 1$, for a.e. $t>0$.

Next, we introduce the time optimal control problem
with control constraints of the rectangular type.
 Arbitrarily fix a sequence of numbers $\{k_i\}_{i=1}^{d}$
 such that $1=k_1\geq k_2\geq\cdots\geq k_d>0$. For each  $M>0$, we define the following  set:
\beqnx
&\mathcal{U}^M\triangleq\left\{u=(u^1,\cdots,u^d):
(0,+\infty)\rightarrow \R^d \q\text{is  measurable};\right.\\
&\left.\q\q\q|u^i(t)|\leq k_i M,\;\;\text{for a.e.}\;\; t\in(0,+\infty) \;\;\text{and for all}\;\;i=1,\cdots,d \right\}.
\eqnx
This set is called
 a control constrained set with the rectangular type. Now, we introduce the following  time optimal control problem:
$$
(TP)^M\q\inf\left\{T>0; y(T;u)=0, u\in\mathcal{U}^M\right\}.
$$
In this problem, the number
 $$t^*(M)\triangleq\inf\left\{T;y(T;u)=0,u\in \mathcal{U}^M \right\}$$
  is called the optimal time; a control $u_*\in\mathcal{U}^M$, such that $y(t^*(M);u_*)=0$, is called a time optimal control (or an optimal control, for simplicity); and
   a control
$u\in\mathcal{U}^M$, such that  $y(T;u)=0$ for some $T>0$, is called an admissible control.

The existence of time optimal controls to $(TP)^M$
has been studied in \cite{SB}. In this paper, we  build up  a necessary and sufficient condition
and provide an algorithm for the optimal time, as well as the optimal control to $(TP)^M$.
To present the first main result, we  introduce a
functional $J^T$, for each $T>0$ and each  $\{k_i\}_{i=1}^{d}$, with  $1=k_1\geq k_2\geq\cdots\geq k_d>0$, by setting
\begin{equation}\label{13}
J^{T}(\varphi_{_T})=\frac{1}{2}\left(\int^T_0\sum_{i=1}^{d}
k_i|\lg b_i,\v(t)\rg|\d t\right)^2+\lg\varphi(0),y_0\rg,\;\;
\;\varphi_{_T}\in \R^m ,
\ee
where $\varphi(\cdot)$ is the solution to Equation (\ref{3}) with $\varphi(T)=\varphi_{_T}$.
It is proved that this functional has a nonzero
minimizer in $\R^m$. Then, the first main result is stated as follows:

\begin{Theorem}\label{Theorem 2}
Let $M>0$. Then, $t^*$ and $u_*$ are the optimal time and the optimal control to $(TP)^M$, respectively, if and only if $t^*>0$ and $u_*\in \mathcal{U}^M$ satisfy that
\bb \label{27}
u^i_*(t)=k_i M\f{\lg b_i,\hat{\varphi}(t)\rg}{|\lg b_i,\hat{\varphi}(t)\rg|}
 \;\;\;\mbox{for a.e.}\;\; t\in(0,t^*) \;\;\mbox{and for all}\;\;i=1,\cdots,d
\ee
and
\bb \label{28}
M=\int_{0}^{t^*}\sum_{j=1}^dk_j|\lg b_j,\hat{\varphi}(t)\rg|\,\d t,
\ee
where \(\hat{\varphi}(\cdot)\) is the solution to Equation (\ref{3}) with \(\varphi(t^*)=\hat{\varphi}_{t^*}\), which is a minimizer of the functional \(J^{t^*}\).
\end{Theorem}

To state the second main result, we define, for each $T>0$ and each $\{k_i\}_{i=1}^{d}$, with  $1=k_1\geq k_2\geq\cdots\geq k_d>0$, a set of controls:
\beqnx
&\mathcal{V}^T\triangleq\left\{v=(v^1,\cdots,v^d)
\in L^\infty(0,T;\R^d);\;\;|v^i(t)|\leq k_i\|v^1\|_{L^\infty(0,T;\R)}
\right.\\
&\q\q\q\q\q\q\q\q\q\q\q\q\q
\text{for a.e.}\;\; t\in(0,T)\;\;
 \text{and for all}\;\;i=1,\cdots,d \big\};
\eqnx
and then introduce the following norm optimal control
problem:
$$(NP)^T\q\inf\left\{\|v^1\|_{L^\infty(0,T;\R)};y(T;v)=0, v\in
\mathcal{V}^T
\right\},$$
where $y(\cdot;v)$ is the solution to Equation (\ref{1}), where the time horizon $\mathbb{R}^+$ is replaced by $(0,T)$,
corresponding to the control $v(\cdot)$.  Write
 $$\widetilde{M}(T)\triangleq\inf\left\{\|v^1\|_
{L^\infty(0,T;\R)};y(T;v)=0,v\in
\mathcal{V}^T\right\}.$$ Then, we  construct
a sequence of numbers $\{t_n\}_{n=0}^{+\infty}$ as follows:
Let $t_0>0$ be arbitrarily given. Let $K\in\mathbb{N}$
be such that
$$K=\min\left\{k\in\mathbb{N};\;\;\widetilde{M}(kt_0)<M\right\}.$$
It is proved that such a $K$ exists. Then we set $a_0=0$, $b_0=Kt_0$. Write $t_1=(a_0+b_0)/2.$
In general, when $t_n=(a_{n-1}+b_{n-1})/2$, $n\geq 1$,
with $a_{n-1}$ and $b_{n-1}$ being given, it is defined
that
\begin{equation*}
[a_{n},b_{n}]=
\begin{cases}
[t_n,b_{n-1}],\;\;\;\text{if}\;\;\widetilde{M}(t_n)>M,\\
[a_{n-1},t_n],\;\;\;\text{if}\;\;\widetilde{M}(t_n)\leq M\\
\end{cases}
\end{equation*}
and $t_{n+1}=(a_n+b_n)/2$. It is proved  that the sequence $\{t_n\}_{n=0}^{+\infty}$  can be determined by solving
a series of problems of calculus of variation $\min_{\varphi_{_T}\in \mathbb{R}^m} J^T(\varphi_{_T})$ with different $T$.
Now, the second main result is presented as follows:

\begin{Theorem}\label{Theorem A}
Suppose that $M>0$.
Let $\{t_n\}_{n=0}^{+\infty}$ be the above-mentioned sequence. Write  $u_*$ for  the optimal control
to  $(TP)^M$. Let  $u_n=(u_n^1,\cdots,u_n^d)$, $n\in\mathbb{N}$, be defined by
$$u^i_n(t)=\left(\int_{0}^{t_n}\sum_{j=1}^d k_j|\lg b_j,
\hat{\v}_n(t)\rg|\,\d t\right)
\frac{k_i\lg b_i,\hat{\v}_n(t)\rg}
{|\lg b_i,\hat{\v}_n(t)\rg|}\;\;\text{for a.e.}\; t\in(0,+\infty),\;\;i=1,\cdots,d,$$
where $\hat{\v}_n(\cdot)$ is the solution to Equation
(\ref{3}) with $\v(t_n)=\hat{\v}_{_{t_n}}$,
which is a minimizer of  the functional $J^{t_n}$.
Then, it holds that
\bb\label{4.1}
t_n\rightarrow t^*(M) \;\;\text{as}\;\; n\rightarrow +\infty
\ee
and for each $i$ with $1\leq i\leq d$,
\bb\label{4.2}
u_n^i\rightarrow u^i_* \;\;\text{in}\;\; L^2(0,t^*(M);\;\R).
\ee
\end{Theorem}

The main idea to prove the above theorems is as follows: We first build up an equivalence theorem of our time optimal control problem and the norm optimal control problem constructed above;
 then make use of the variational characterization of the norm optimal control and the equivalence theorem to show the above two theorems.  The aforementioned  equivalence theorem is motivated by the analogous equivalence results  established
for heat equations with control constraints of the ball type in \cite{WZ} (see also \cite{WX}). However, the time optimal control problems with control constraints of the rectangular type  differ from
those with control constraints of the ball type, from different points of view (see for instance \cite{LW}).   The equivalence theorem, as well as the structure of the norm optimal
 control problems in this paper seems to be new.

There have been  a lots of literatures on time optimal control problems of differential equations (see, for instance, \cite{EV}, \cite{FA}, \cite{LPW},  \cite{LM},  \cite{LW}, \cite{MSE}, \cite{PB}, \cite{SO}, \cite{SB}, \cite{W} ). Recently, the  semi-smooth Newton methods to analyze numerically the time optimal controls with  constraints of the cubic type
for some ordinary differential equations have been introduced
in \cite{KU}.

To  our best
knowledge, the necessary and sufficient
condition and the algorithm for the optimal time, as well as optimal control, provided in this paper,
have not been studied.
The equality (\ref{28}) provides a formula
for the optimal time to $(TP)^M$.

The rest of this paper is organized as follows: Section 2 provides some results related to the norm optimal control problem $(NP)^T$. Section 3 establishes an equivalence theorem between the norm and the time optimal controls.  Section 4 presents the proof of the main theorems.

\section{Some Properties about $(NP)^T$}

$\;\;\;\;$We first present the following properties for the functional $J^T$ which is defined by (\ref{13}).
\begin{Lemma}\label{Lemma 3}
For each $T>0$, $J^T$ is continuous, convex, and
coercive in $\R^m$. Moreover, each minimizer of $J^T$ is nonzero.
\end{Lemma}
\begin{proof} We first show the existence of minimizers for $J^T$. The proof of the continuity and convexity of $J^T$ follows from the standard argument (see, for instance,  \cite{Z}, \cite{ZZ}). Now, we  show the coercivity of $J^T$. For this purpose,   we  set, for each $\v_{_T}\in \R^m$,
\bb\label{14}
\|\varphi_{_T}\|_*\triangleq\int^T_0\sum_{i=1}^{d}
k_i|\lg b_i,\v(t)\rg|\d t,\ee
where $\v(\cdot)$ is the solution to Equation
(\ref{3}) with $\v(T)=\v_{_T}$. Because of $(H.1)$ and $(H.2)$,  $\|\cdot\|_*$ is a norm in $\R^m$. By the equivalence of norms in $\R^m$, there exists a constant $\Lambda>0$ such that
$$
\|\varphi_{_T}\|\leq \Lambda\|\v_{_T}\|_*.
$$
This, together with the definition of $J^T$ and (\ref{14}), leads to
\begin{eqnarray*}
J^T(\varphi_{_T})&\geq&\frac{1}{2\Lambda^2}\|\v_{_T}\|^2-
\|y_0\|\|\v(0)\| .
\end{eqnarray*}
Thus
$$\lim_{\|\v_{_T}\|\rightarrow+\infty}J^T(\v_{_T})=+\infty.$$
Hence,  $J^T(\cdot)$ is coercive in $\R^m$. Therefore,  $J^T$ has  minimizers
in $\R^m$.

We next show that any minimizer of $J^T$ is nonzero. For this purpose, we set, for each $\alpha>0$,  $\varphi^\alpha_{_T}\triangleq-\alpha \Psi(T,0)y_0$, where $\Psi(\cdot,\cdot)$ is the fundamental solution associated to Equation (\ref{3}).
Write $\varphi^\alpha(\cdot)$ for the solution to Equation (\ref{3}) with $\v(T)=\varphi^\alpha_{_T}$. Then
\begin{eqnarray*}
J^T(\varphi^\alpha_{_T})
&=&\frac{\alpha^2}{2}\left(\int_0^T\sum_{i=1}^{d}k_i|\lg b_i,\Psi(t, 0)y_0\rg|\d t\right)^2
-\alpha\|y_0\|^2\\
&\leq&\frac{\alpha^2T^2}{2}\|\Psi(\cdot,0)\|^2_{L^{\infty}
(0,T;\R^{m\t m})}\left(\sum_{i=1}^{d}k_i\|b_i\|\right)^2\|y_0\|^2
-\alpha\|y_0\|^2.
\end{eqnarray*}
 This, together with the assumption that $y_0\neq 0$, implies that $J^T(\varphi^\alpha_{_T})<0$ whenever $\alpha>0$ small enough. Therefore, each minimizer is nonzero. This completes the proof.
\end{proof}
\begin{Remark}\label{remarkj}
In general, the functional $J^T(\cdot)$ defined by
(\ref{13}) is not strictly convex in $\mathbb{R}^m$.
Here we give an example to explain it. Now, assume
that $m=2$, $d=1$, $T=\pi/4$, $k_1=1$ and
\[A(\cdot)\equiv A=
\left(\begin{array}{cc}
0 &-1\\
1 & 0
\end{array}\right),\quad
 b_1=
\left(\begin{array}{c}
1\\
0
\end{array}\right).
\]
Then $(A,b_1)$ satisfies Conti's condition. Next, set
\[
\varphi^1_{_T}=\left(\begin{array}{c}
1\\
0
\end{array}\right),\quad
\varphi^2_{_T}=\left(\begin{array}{c}
1+\frac{\sqrt{2}}{2}\\
1+\frac{\sqrt{2}}{2}
\end{array}\right).
\]
Let $\varphi^i(\cdot)$ be the solution to Equation (\ref{3}) with $\varphi(T)=
\varphi^i_{_T}$, $i=1,2$, respectively. Then one can readily
check that for each $\lambda\in(0,1)$, $J^T\big(\lambda\varphi^1_{_T}+(1-\lambda)\varphi^2_{_T}\big)
=\lambda J^T(\varphi^1_{_T})+(1-\lambda)J^T(\varphi^2_{_T})$.
Thus $J^T(\cdot)$ is not strictly convex in $\mathbb{R}^2$.
\end{Remark}

\begin{Lemma}\label{Lemma 3.1}
 Suppose that $T>0$. Let $\hat{\v}_{_T}$ be a minimizer of $J^T$. Write  $\hat{\v}(\cdot)$ for the solution to Equation (\ref{3}) with
$\v(T)=\hat{\v}_{_T}$.
Then, the control $\bar{u}=(\bar{u}^1,\cdots,\bar{u}^d)$, where
\bb\label{15}
\bar{u}^i(t)=\left(\int_0^T\sum_{j=1}^{d}k_j|\lg b_j,\hat{\v}(t)\rg|\d t\right)
\frac{k_i\lg b_i,\hat{\v}(t)\rg}{|\lg b_i,\hat{\v}(t)\rg|},
\;\;t\in(0,T),\; \;i=1,\cdots d,
\ee
 is   optimal
to  $(NP)^T$. Consequently,
\bb\label{16}
\widetilde{M}(T)=\int_0^T\sum_{j=1}^{d}k_j|\lg b_j,\hat{\v}(t)\rg|\d t.
\ee
\end{Lemma}
\begin{proof}
According to Lemma~\ref{Lemma 3}, it holds that $\hat{\v}_{_T}\neq0$.  This,
 along with $(H.1)$ and $(H.2)$, indicates that
  $\lg b_i,\hat{\v}(t)\rg\neq0$ for a.e. $t\in(0,T)$ and for all $i=1,\cdots,d$.

Next, we derive the Euler-Lagrange equation
of the functional $J^T$ associated with $\hat{\v}_{_T}$.
For each $\v_{_T}\in \R^m$, let $\v(\cdot)$ be the solution to Equation (\ref{3}) with $\v(T)=\v_{_T}$.
Then
\begin{equation*}
\begin{split}
&\lim_{h\rightarrow0}\frac{1}{h}
\left[J^T(\hat{\v}_{_T}+h\v_{_T})
-J^T(\hat{\v}_{_T})\right]
=\lg\v(0),y_0\rg+\\
&\lim_{h\rightarrow 0}
\frac{1}{2h}\left[\left(\int^T_0\sum_{i=1}^{d}
k_i|\lg b_i,\hat{\v}(t)+h\v(t)\rg|\d t\right)^2-\left(\int^T_0\sum_{i=1}^{d}
k_i|\lg b_i,\hat{\v}(t)\rg|\d t\right)^2\right]\\
&=\left(\int^T_0\sum_{i=1}^{d}
k_i|\lg b_i,\hat{\v}(t)\rg|\d t\right)\int^T_0\sum_{i=1}^{d}
\frac{k_i\lg b_i,\hat{\v}(t)\rg}{|\lg b_i,\hat{\v}(t)\rg|}\lg b_i,\v(t)\rg\;\d t+\lg\v(0),y_0\rg.
\end{split}
\end{equation*}
Hence,  for each $\v_{_T}\in\R^m$, it stands that
\bb\label{17}
\left(\int^T_0\sum_{i=1}^{d}
k_i|\lg b_i,\hat{\v}(t)\rg|\d t\right)\int^T_0\sum_{i=1}^{d}
\frac{k_i\lg b_i,\hat{\v}(t)\rg}{|\lg b_i,\hat{\v}(t)\rg|}\lg b_i,\v(t)\rg\;\d t+\lg\v(0),y_0\rg=0.
\ee

The remainder is to show that  $\bar{u}$ is an optimal control
 to $(NP)^T$. For this purpose, we first observe that  for each   $v\in\mathcal{V}^T$,  $y(T;v)=0$ if and only if
\bb\label{18}
\int_{0}^{T}\sum_{i=1}^{d}v^i(t)\lg b_i,\v(t)\rg\d t
+\lg\v(0),y_0\rg=0,\;\;\text{for each}\;\;\v_{_T}\in\R^m.
\ee
On the other hand, it follows from  (\ref{15}) that
$\bar{u}\in\mathcal{V}^T$, while it follows from  (\ref{15}), (\ref{17}) and (\ref{18})  that $y(T;\bar{u})=0$.
Now, by taking $\v_{_T}=\hat{\v}_{_T}$
 in both (\ref{17}) and (\ref{18}) respectively, and then using  (\ref{15}), we  deduce that
\bb\label{21}
\int_{0}^{T}\sum_{i=1}^{d}\bar{u}^i(t)\lg b_i,\hat{\v}(t)\rg\d t=\int_{0}^{T}\sum_{i=1}^{d}v^i(t)\lg b_i,\hat{\v}(t)\rg\d t.
\ee
By making use of (\ref{15}) again, we see that
\bb\label{22}
\int_{0}^{T}\sum_{i=1}^{d}\bar{u}^i(t)\lg b_i,\hat{\v}(t)\rg\d t=\left(\int_0^T\sum_{i=1}^{d}k_i|\lg b_i,\hat{\v}(t)\rg|\d t\right)^2.
\ee
 Since $v\in\mathcal{V}^T$, it stands that
\begin{equation*}
\left|\int_{0}^{T}\sum_{i=1}^{d}v^i(t)\lg b_i,\hat{\v}(t)\rg\d t\right|\leq
\|v^1\|_{L^{\infty}(0,T;\R)}
\left(\int_0^T\sum_{i=1}^{d}k_i|\lg b_i,\hat{\v}(t)\rg|\d t\right).
\end{equation*}
This, combined  with (\ref{21}) and (\ref{22}), yields that
\begin{equation*}
\int_0^T\sum_{i=1}^{d}k_i|\lg b_i,\hat{\v}(t)\rg|\d t\leq
\|v^1\|_{L^{\infty}(0,T;\R)},
\end{equation*}
from which and (\ref{15}), it follows that
$$\|\bar{u}^1\|_{L^{\infty}(0,T;\R)}
\leq \|v^1\|_{L^{\infty}(0,T;\R)}.
$$
This completes the proof.\\
\end{proof}

\section{Equivalence of Time and Norm Optimal Controls }

$\;\;\;\;\;$The main purpose of this section is to show the following equivalence theorem:
\begin{Theorem}\label{Theorem 1}
For each \(T>0\), the norm optimal control  to $(NP)^T$, when  is extended to $(0,+\infty)$ by taking zero value on $(T,+\infty)$,  is  the time optimal control to \((TP)^{\widetilde{M}(T)}\). On the other hand, for each \(M>0\), the time optimal control to $(TP)^M$, when is restricted over $(0,t^*(M))$,  is  the norm optimal control to $(NP)^{t^*(M)}$ .
\end{Theorem}
\noindent We start with introducing  three lemmas as follows:

\begin{Lemma}\label{Lemma 1}
For each $M>0$, $(TP)^M$ has a unique optimal control over $(0,t^*(M))$. Furthermore, it has the
strong bang-bang property: any optimal control $u_*=(u_*^1,\cdots, u_*^d)$ satisfies
that  $|u_*^i(t)|=k_i M$, for a.e. $t\in(0,t^*(M))$ and for all $i=1,\cdots,d$.
\end{Lemma}
\begin{proof}
Since $(H.2)$ stands, there exists a control $u\in\mathcal{U}^M$ such that $y(T;u)=0$ for some $T>0$ (see \cite{CON}, \cite{SB}).
By a standard  argument (see, for instance, \cite{LM}), the existence
of  time optimal controls to  $(TP)^M$ follows immediately.

Next, let $u_*$ be an optimal control to  $(TP)^M$.
By the Pontryagin maximum principle (see, for instance, \cite{PB}), there
exists a nonzero solution $\v(\cdot)$ to Equation (\ref{3}) such that
\bb\label{4}\sum_{i=1}^d\max_{|v^i|\leq k_iM}\lg b_i,\v(t)\rg v^i=
\sum_{i=1}^d\lg b_i,\v(t)\rg u_*^i \;\;\text{for a.e.}\;\; t\in(0,t^*(M)).
\ee
Because of  $(H.1)$ and $(H.2)$, it holds that  $\lg b_i,\v(t)\rg\neq0$ for
a.e. $t\in(0,t^*(M))$ and for  all $i=1,\cdots,d$. This, together with (\ref{4}), yields that
$$u_*^i(t)=k_iM\frac{\lg b_i,\v(t)\rg}{|\lg b_i,\v(t)\rg|}
\;\;\text{for a.e.}\;\;t\in(0,t^*(M))\;\;\text{and for all}\;\;i=1,\cdots,d.$$
Hence, the desired strong bang-bang property follows
immediately. Finally, by the strong bang-bang property, the uniqueness of the time optimal control over $(0,t^*(M))$ follows from the standard argument (see, for instance, \cite{FA}).
\end{proof}

\begin{Lemma}\label{lemma 2}
Let $T>0$. For each $\tau\in[0,T)$ and each $z_0\in\R^{m}$, there exists a control $u^1\in L^{\infty}(\tau,T;\R)$
such that the solution $y(\cdot;u^1)$ to the following equation:
\bb\label{5}
y'(t)+A(t)y(t)=b_1u^1(t),\;\;
y(\tau)=z_0,
\ee
verifies $y(T;u^1)=0$.
Moreover, the control $u^1$ satisfies the following
estimate:
$$\|u^1\|_{L^{\infty}(\tau,T;\R)}
\leq C\|z_0\|,$$
where $C$ is a positive constant independent of $z_0$.
\end{Lemma}
\begin{proof}
Since $(H.1)$ stands and $(A(\cdot),b_1)$ satisfies  Conti's
condition (see (H.2)), the system (\ref{5}) holds the unique
continuation property on $(\tau, T)$. Then, applying the
Theorem 5 in Chapter 3 in \cite{SO}, we get
that the controllability Gramian
$W(\tau,T)$ is positive definite, where
$$W(\tau,T)=\int_\tau^T\Phi(T,s)b_1b_1^*\Phi^*(T,s)ds.$$
Here, $\Phi(\cdot,\cdot)$ is the fundamental solution associated to $A(\cdot)$.  Next, set
\bb\label{+11}
u^1(t)=-b_1^*\Phi^*(T,t)W(\tau,T)^{-1}\Phi(T,\tau)z_0,
\;\;t\in[\tau,T).
\ee
It can be easily checked that $y(T;u^1)=0$.
By (\ref{+11}), it holds that
$$\|u^1\|_{L^\infty(\tau,T;\R)}\leq\|b_1^*
\Phi^*(T,\cdot)\|_{L^\infty(\tau,T;\R^m)}\|W
(\tau,T)^{-1}\|_{\R^{m\times m}}\|
\Phi(T,\tau)\|_{\R^{m\times m}}\|z_0\|.$$
 Hence, there exists a constant $C>0$ (independent of $z_0$)  such that
$$\|u^1\|_{L^\infty(0,T;\R)}\leq C\|z_0\|.$$
This completes the proof.
\end{proof}

\noindent Next lemma concerns some properties of the
map $M\rightarrow t^*(M)$.
\begin{Lemma}\label{Lemma 3.2}
The optimal time function $t^*(\cdot)$ is strictly monotonically decreasing and continuous. In addition, it holds that
$\lim\limits_{M\to +\infty}t^*(M)=0$ and  $\lim\limits_{M\to 0^+}t^*(M)=+\infty$.
\end{Lemma}
\begin{proof}
We carry out the proof by five steps as follows:\\

\noindent {\it  Step 1: The function \(t^*(\cdot)\) is strictly monotonically decreasing.}

Let \(M_1> M_2>0\). It suffices to show that $t^*(M_1)<t^*(M_2)$. To this end, let  $u_2$ be the optimal control to   \((TP)^{M_2}\). Clearly,  \(u_2\) is admissible for \((TP)^{M_1}\).
By the optimality of \(t^*(M_1)\) to \((TP)^{M_1}\), it is clear that \(t^*(M_1)\leq t^*(M_2)\). Next, suppose by contradiction that \(t^*(M_1)=t^*(M_2)\). Then,  it would  hold that
\[
y(t^*(M_1);u_2)=y(t^*(M_2);u_2)=0.
\]
Hence, \(u_2\) is also the optimal control to
 \((TP)^{M_1}\).
By Lemma~\ref{Lemma 1}, we find that
$$
k_iM_1=|u^i_2(t)|=k_iM_2\;\;\text{for a.e.}\;\; t\in (0,t^*(M_1))
\;\;\text{and for all}\;\;i=1,\cdots,d.
$$
This leads to a contradiction, since \(M_2<M_1\). Therefore,
it holds that \(t^*(M_1)< t^*(M_2)\).\\

\noindent {\it Step 2: The function \(t^*(\cdot)\) is continuous from right,
that is, $\lim\limits_{M_n\searrow M}t^*(M_n)=t^*(M)$.}

Let \(M_1> M_2> \cdots> M_n>\cdots>M>0\) and \(\lim\limits_{n\to +\infty}M_n=M\). By Step 1, it holds that
$$t^*(M_1)< t^*(M_2)<\cdots< t^*(M_n)<\cdots<t^*(M)\;\;\text{ and}\;\;\lim_{n\to +\infty} t^*(M_n)\leq t^*(M).
$$
We claim that $\lim\limits_{n\to +\infty} t^*(M_n)=t^*(M)$. Seeking a contradiction, we suppose that
\begin{eqnarray*}
\lim\limits_{n\to+\infty} t^*(M_n)= t^*(M)-\delta\;\;\mbox{for some}\;\; \delta>0.
\end{eqnarray*}
Clearly, the optimal controls $u_n$ to
$(TP)^{M_n}$, $n\in\mathbb{N}$,
satisfy that
\begin{equation*}
\|u^i_n\|_{L^{\infty}(\mathbb{R}^+;\;\R)}\leq k_iM_n< k_iM_1,\;\; \text{for all}\;\;i=1,\cdots,d
\end{equation*}
and
$$
y(t^*(M_n); u_n)=0.
$$
Thus, on a subsequence, $u_n\rightarrow\widetilde{u}$ weakly star in $L^\infty(\mathbb{R}^+; \mathbb{R}^d)$. Furthermore, one can easily derive from the above observations that
$\widetilde{u}\in \mathcal{U}^M$ and
$y(t^*(M)-\delta; \widetilde{u})=0$. These contradict
with the optimality of $t^*(M)$ to $(TP)^M$.\\

\noindent {\it Step 3: The function \(t^*(\cdot)\) is continuous from left, that is, $\lim\limits_{M_n\nearrow M}t^*(M_n)=t^*(M)$.}

Let \(0<M_1< M_2< \cdots< M_n<\cdots<M\) and \(\lim\limits_{n\to\infty}M_n=M\). It is clear that
$$t^*(M_1)> t^*(M_2)>\cdots> t^*(M_n)>\cdots>t^*(M)\;\;\text{ and}\;\;\lim_{n\to+\infty} t^*(M_n)\geq t^*(M).$$

Seeking a contradiction, suppose $\lim\limits_{n\to+\infty} t^*(M_n)>t^*(M)$. Then there would  exist a \(\delta>0\) such that \[\lim_{n\to+\infty} t^*(M_n)= t^*(M)+\delta.\]
 Clearly,
\bb\label{+3}
t^*(M_n)>t^*(M)+\delta\;\;\text{for all}\;\;n\in\mathbb{N}.
\ee
Let $u_*$ be the optimal control to
$(TP)^{M}$.
Set \(\delta_n=\frac{M_n}{M}\), \(z_n(\cdot)=\delta_n y(\cdot;u^*)\), $n\in\mathbb{N}$.
It is clear that
\bb\label{+9}
\begin{cases}
z_n'(t)+A(t)z_n(t)=\sum_{i=1}^{d}\delta_n b_i u_*^i(t),\;\;t\in(0,t^*(M)),\\
z_n(0)=\delta_n y_0,\q z_n(t^*(M))=0.
\end{cases}
\ee
According to  Lemma~\ref{lemma 2}, there exist a constant $C>0$ independent of $n$ and  a control \(f^1_n\) with
$$
\|f^1_n\|_{L^{\infty}(t^*(M),t^*(M)+
\delta;\;\R)}\leq C\cdot(1-\delta_n)\| y_0\|,
$$
such that
\[
\phi_n(t^*(M)+\delta)=0,
\]
where $\phi_n(\cdot)$ solves the equation:
\bb\label{+10}
\begin{cases}
\phi'(t)+A(t)\phi(t)=b_1f^1_n(t)\chi_{_{(t^*(M), t^*(M)+\delta)}}(t),\;\;t\in(0,t^*(M)+
\delta),\\
\phi(0)=(1-\delta_n)y_0.
\end{cases}
\ee

Now, we construct, for each $n\in\mathbb{N}$,  a control  $g_n=(g_n^1,\cdots,g_n^d)$, by setting
\begin{equation}\label{+4}
\begin{cases}
g^1_n=\delta_nu_*^1\chi_{_{(0,t^*(M))}}+f^1_n\chi
_{_{(t^*(M),t^*(M)+\delta)}},\\
g^i_n=\delta_nu_*^i\chi_{(0,t^*(M))},\;\; i=2,\cdots,d.
\end{cases}
\end{equation}
Since $\delta_n\nearrow 1$, there exists a positive integer
$N_1$ such that
\[
C\cdot(1-\delta_n)\|y_0\|\leq M_1< M_n\;\;\text{for all}\;\;
 n\geq N_1.
\]
This, along with (\ref{+4}), leads to that when $n\geq N_1$
\[
\|g^i_n\|_{L^{\infty}(\R^+;\;\R)}\leq k_i M_n,
\;\;\text{for all}\;\;i=1,\cdots,d.
\]
Finally, set \(w_n=z_n+\phi_n\), $n\geq N_1$.
It follows at once from (\ref{+9}) and (\ref{+10}) that
\[
\begin{cases}
w_n'+A(t)w_n(t)=\sum_{i=1}^{d}b_ig_n^i(t), \q t\in(0,t^*(M)+\delta),\\
w_n(0)=y_0,\q w_n(t^*(M)+\delta)=0.
\end{cases}
\]
Thus, \(g_n\) is admissible to \((TP)^{M_n}\) for each $n$ with $n\geq N_1$. Consequently, $t^*(M_n)\leq
t^*(M)+\delta$ whenever $n\geq N_1$.
This, together with (\ref{+3}), leads to a contradiction.\\

\noindent {\it Step 4: It holds that  \(\lim\limits_{M\to 0^+}t^*(M)=+\infty\).}

Seeking a contradiction,  suppose that there did exist a sequence $\{M_n\}_{n\geq 1}$, with $M_1> M_2> \cdots> M_n>\cdots>0$ and $\lim\limits\limits_{n\to \infty}M_n=0$, such that $\lim_{n\to\infty} t^*(M_n)= T<+\infty$.
Then, the optimal controls  $u_n$  to \((TP)^{M_n}\), $n\in\mathbb{N}$, satisfy that on a subsequence,
$y(\cdot;u_n)\to y(\cdot;0)$ in $C([0,T];\R^d)$, which leads to a contradiction, since  $y_0\neq0$.\\

\noindent {\it  Step 5: \(\lim\limits_{M\to + \infty}t^*(M)=0\)}

Seeking a contradiction, suppose that  there existed a
   $T>0$ and a sequence \(\{M_n\}_{n\geq 1}\), with
\(0<M_1< M_2<\cdots< M_n<\cdots \;\textrm{and }\lim\limits_{n\to\infty}M_n= +\infty\), such that $\lim\limits_{n\to+\infty}t^*(M_n)=T$.\\
Let \(\delta>0\) such that \(T-\delta>0\). Then, by
 Lemma~\ref{lemma 2}, there would exist  a control \(u^1_\delta\)
with
$$
\|u^1_\delta\|_{L^{\infty}(0, T-\delta;\;\R)}\leq C\|y_0\|,
$$
such that
\bb\label{12}
y(T-\delta;(u^1_\delta,\underbrace{0,\cdots,0}_{d-1}))=0.
\ee
Since \(\lim\limits_{n\to +\infty}M_n=+\infty\), it holds that $C\|y_0\|\le M_n$
 for $n$ large enough. This, together with (\ref{12}), leads
to a contradiction to the optimality of \(t^*(M_n)\) to
\((TP)^{M_n}\).\\

In summary, we conclude that all statements in this lemma stand.\\
\end{proof}

\noindent\textbf{Proof of Theorem~\ref{Theorem 1}.}$\;\;$
We begin with proving the identity
\begin{equation}\label{24}
 T=t^*(\widetilde{M}(T))\;\;\text{for each}\;\; T>0.                                        \end{equation}
From the definition of $\widetilde{M}(T)$ and the optimality of \;\(t^*(\widetilde{M}(T))\) to \((TP)^{\widetilde{M}(T)}\), we can deduce that for each
$T>0$, $t^*(\widetilde{M}(T))\leq T$.
Thus, it suffices to show that the inequality $t^*(\widetilde{M}(T))<T$ does not stand for each $T>0$.
Suppose by contradiction that \(t^*(\widetilde{M}(T))< T\)
for some $T>0$. Then, by making use of Lemma~\ref{Lemma 3.2},
we could find a positive number $M_1$, with \(M_1<\widetilde{M}(T)\), such that \(t^*(M_1)= T\). It follows from Lemma~\ref{Lemma 1} that
$(TP)^{M_1}$ has a unique optimal control $u_*$ verifying
\begin{equation*}
|u_*^i(t)|=k_i M_1 \;\;\text{for a.e.}\;\;
t\in(0,T)\;\;\text{and for all}\;\;i=1,\cdots,d.
\end{equation*}
Thus, $u_*\in \mathcal{V}^T$ and $y(T;u_*)=0$.
This contradicts with  the optimality
of $\widetilde{M}(T)$ to
$(NP)^T$. Therefore, the equality (\ref{24}) stands.

Now,  any optimal control  $f_*$  to $(NP)^T$ satisfies that
$|f_*^i(t)|\leq k_i\widetilde{M}(T)$ for a.e. $t\in(0,T)$ and all $i=1,\cdots,d$, and $ y(T;f_*)=0$.
These, along with (\ref{24}), imply that
$
f_*\in\mathcal{U}^{\widetilde{M}(T)}$ and $ y(t^*(\widetilde{M}(T));f_*)=0$.
Hence, \(f_*\) is the  optimal control to \((TP)^{\widetilde{M}(T)}\).

On the other hand,  it follows from (\ref{24}) and  the strict monotonicity of
the function $t^*(\cdot)$ that
\bb\label{+7}
\widetilde{M}(t^*(M))=M,\;\;\text{for each}\;\;M>0.
\ee
Thus, the optimal control $u_*$  to $(TP)^M$
is the optimal control to \((TP)^{\widetilde{M}(t^*(M))}\). Then, by the optimality of $u_*$ and by Lemma~\ref{Lemma 1}, we see
that $\|u_*^1\|_{L^\infty(0,t^*(M);\;\R)}=\widetilde{M}(t^*(M))$,
$u_*\in\mathcal{V}^{t^*(M)}$ and  $y(t^*(M);u_*)=0$. Hence,
$u_*$ is an optimal control to  $(NP)^{t^*(M)}$.
This completes the proof of Theorem~\ref{Theorem 1}.\\

We end this section with the following two consequences.

\begin{Corollary}\label{Theorem*}
For each \(T>0\), $(NP)^T$ has a unique optimal control \(f_*=(f^1_*,\cdots, f^d_*)\). It is given by
$$
f_*^i(t)=\left(\int_0^T\sum_{j=1}^{d}
k_j|\lg b_j,\hat{\v}(t)\rg|\d t\right)
\frac{k_i\lg b_i,\hat{\v}(t)\rg}{|\lg b_i,\hat{\v}(t)\rg|}, \;\text{for a.e. }\;t\in(0, T)\;
\text{and for all}\;\;i=1,\cdots,d,
$$
where $\hat{\v}(\cdot)$ is the solution to Equation (\ref{3}) with
$\v(T)=\hat{\v}_{_T}$, which is a minimizer of the functional $J^T$. Consequently,
$$\widetilde{M}(T)=\int_0^T\sum_{j=1}^{d}
k_j|\lg b_j,\hat{\v}(t)\rg|\d t.$$
\end{Corollary}
\begin{proof}
It  suffices to show the uniqueness, because of  Lemma~\ref{Lemma 3.1}.
For this purpose, we   suppose, by contradiction,  that \(g_*\) were an optimal control different from $f_*$. Then, according to  Theorem~\ref{Theorem 1},
 both  $f_*$ and $g_*$ were optimal controls to \((TP)^{\widetilde{M}(T)}\). By Lemma~\ref{Lemma 1}, as well as (\ref{24}), they are the same over $(0,T)$, which leads to a contradiction.
  This completes the proof.
\end{proof}

\begin{Corollary}\label{Prop 2}
The functions \(\widetilde{M}(\cdot)\) and  \(t^*(\cdot)\) are  inverse one of  each other. Consequently, $\widetilde{M}(\cdot)$ is strictly monotonically decreasing and continuous. In addition,
it holds that
$$
\lim_{T\rightarrow0^+}\widetilde{M}(T)=+\infty
\;\;\text{and}\;\; \lim_{T\rightarrow+\infty}\widetilde{M}(T)=0.
$$
\end{Corollary}
\begin{proof}
According to Theorem~\ref{Theorem 1},
it follows that
$$t^*\circ\widetilde{M}(T)=t^*(\widetilde{M}(T))=T,
\;\;\text{for each}\;\;T>0$$
and
$$\widetilde{M}\circ t^*(M)=\widetilde{M}(t^*(M))=M,
\;\;\text{for each}\;\;M>0.$$
Hence, $\widetilde{M}(\cdot)$ is the inverse function
of $t^*(\cdot)$. The remainders follow at once
from Lemma~\ref{Lemma 3.2}. This completes the proof.
\end{proof}

\section{Proof of Theorem~\ref{Theorem 2} and \ref{Theorem A}}

\noindent\textbf{Proof of Theorem~\ref{Theorem 2}.}
\;Suppose that
$t^*$ and
$u_*$ are the optimal time and the optimal control to
$(TP)^M$ respectively. It is clear that
$t^*=t^*(M)$. This, combining with
Theorem~\ref{Theorem 1}, yields that
$u_*$ is an optimal control to  $(NP)^{t^*}$.
According to Corollary~\ref{Theorem*}, each $u^i_*(\cdot)$, with $i\in \{1,\cdots,d\}$, satisfies
\bb\label{+8}
u^i_*(t)=\left(\int_0^{t^*}\sum_{j=1}^{d}k_j|\lg b_j,\hat{\v}(t)\rg|\d t\right)
\frac{k_i\lg b_i,\hat{\v}(t)\rg}{|\lg b_i,\hat{\v}(t)\rg|}
\;\;\text{for a.e.}\;t\in(0,t^*),
\ee
where \(\hat{\varphi}(\cdot)\) is the solution to Equation (\ref{3}) with \(\varphi(t^*)=\hat{\varphi}_{t^*}\),  which is a minimizer of the functional \(J^{t^*}\).
Consequently,
$$
\widetilde{M}(t^*)=\int_0^{t^*}\sum_{j=1}^{d}k_j|\lg b_j,\hat{\v}(t)\rg|\,\mathrm{d}t.
$$
This, along with the fact $t^*=t^*(M)$ and
the identity (\ref{+7}),
leads to (\ref{28}).
The equality (\ref{27}) follows immediately from (\ref{28}) and (\ref{+8}).

Conversely, suppose that $t^*>0$ and $u_*\in \mathcal{U}^M$
satisfy the equality (\ref{27}) and (\ref{28}).
We are going to show that they are the optimal time and the optimal control to
 $(TP)^M$ respectively. For this purpose, we apply Corollary~\ref{Theorem*} to obtain that
$u_*$ is the unique optimal control to
 $(NP)^{t^*}$. It is clear that
$\widetilde{M}(t^*)=M$. By the strict monotonicity of $\widetilde{M}(\cdot)$ (see Corollary~\ref{Prop 2}) and the equality (\ref{+7}),
it holds that $t^*=t^*(M)$. According to
Theorem~\ref{Theorem 1},  $u_*$ is the
 optimal control to  $(TP)^{\widetilde{M}(t^*(M))}$. This, along with
(\ref{+7}), yields that $u_*$ is the
 optimal control to  $(TP)^M$ and
completes the proof of Theorem~\ref{Theorem 2}.\\

\begin{Remark} Before giving the proof of Theorem~\ref{Theorem A},
we explain the well-posedness of the
 sequence $\{t_n\}_{n=0}^{+\infty}$ built up in Section 1.
In fact, for each  $T>0$, we can determine the value
$\widetilde{M}(T)$ by solving the minimization problem: $\min_{\varphi_T\in \mathbb{R}^m}J^T(\varphi_T)$
(see Corollary~\ref{Theorem*}).  Since the map
$T\rightarrow\widetilde{M}(T)$ is strictly
monotonically decreasing and $\widetilde{M}(T)$
tends to $0$ as T goes to $+\infty$ (see
Corollary~\ref{Prop 2}), $K$ can be confirmed
by solving a finite number of  minimizers of
 functionals $J^T$, corresponding to
$T=lt_0$, $l=1,2,\cdots,K$. On the other hand, for
each $n\in\mathbb{N}$,
$t_{n+1}$ is determined by solving the same minimization problem with
$T={t_n}$.
Hence,  the sequence $\{t_{n}\}
_{n=0}^{+\infty}$ can be determined by
solving a series of minimizers of
functionals $J^T$ with $T=lt_0$,
 $l=1,2,\cdots,K$, and $T=t_n$, $n=1,2,\cdots$.
 \end{Remark}

\noindent\textbf{Proof of Theorem~\ref{Theorem A}.}
\;We start with proving (\ref{4.1}). From the
structure of the sequence $\{t_n\}_{n=0}^{+\infty}$,
it is clear that $t_n\in[a_n,b_n]\subseteq[a_{n-1},b_{n-1}]$
 and $b_n-a_n=(b_{n-1}-a_{n-1})/2$. Hence, it stands
 that
 \bb\label{4.3}
\lim_{n\rightarrow +\infty}t_n=\lim_{n\rightarrow +\infty}a_n
=\lim_{n\rightarrow +\infty}b_n.
\ee
Since the function $\widetilde{M}(\cdot)$ is
continuous (see Corollary~\ref{Prop 2})
and $\widetilde{M}(a_n)>M\geq \widetilde{M}(b_n)$
(which also follows from the structure of $\{t_n\}_{n=0}^{+\infty}$), we see that
$$M=\widetilde{M}(\lim_{n\rightarrow +\infty}t_n).$$
This, together with (\ref{+7}) and the strict
monotonicity of the function $\widetilde{M}(\cdot)$,
 leads to the desired convergence (\ref{4.1}).

 Next, we  claim that for each $i=1,\cdots,d$,
 \bb\label{4.4}
 u^i_n\rightarrow u_*^i \;\;\text{weakly star in}\;\; L^\infty(0,t^*(M);\,\R).
 \ee
 In fact, for each $n\in\mathbb{N}$,
 it follows from Corollary~\ref{Theorem*}
 that $u_n$, when is restricted over $(0,t_n)$, is
 the unique optimal control to
 $(NP)^{t_n}$. Consequently, $y(t_n;u_n)=0$. We arbitrarily
 take a subsequence of $\{u_n\}$, denoted
 by $\{u_{n_{k}}\}$. Clearly, there exists a
 subsequence $\{u_{n_{j}}\}$ of
 $\{u_{n_{k}}\}$ such that
 for each $i$ with $1\leq i\leq d$,
 \bb\label{4.5}
 u^i_{n_{j}}\rightarrow \tilde{u}^i \;\;
 \text{weakly star in}\;\; L^\infty(0,t^*(M);\,\R).
 \ee
Moreover, one can derive from the above facts that $y(t^*(M);\tilde{u})=0$.
On the other hand, by (\ref{4.1}) and (\ref{+7}),
 it follows that for each $i=1,\cdots,d$,
 $$\|\tilde{u}^i\|_{L^\infty(0,t^*(M);\R)}
 \leq \liminf_{j\rightarrow+\infty}\|u_{n_{j}}^i\|
 _{L^\infty(0,t^*(M);\R)}=\liminf_{j\rightarrow+\infty}
 k_i\widetilde{M}(t_{n_{j}})=k_i\widetilde{M}
 (t^*(M))=k_iM.$$
 Hence, $\tilde{u}$ is an optimal control to
 $(TP)^M$. By the uniqueness of the
 optimal control to $(TP)^M$,
 it holds that $\tilde{u}=u_*$ over $(0, t^*(M))$. Therefore, (\ref{4.4})
 follows from (\ref{4.5}).

 Now, we verify the  convergence (\ref{4.2}).
 By (\ref{4.4}),
  we find that
 for each $i$ with $1\leq i\leq d$,
 \bb\label{4.6}
 u_n^i\rightarrow u_*^i \;\;\text{weakly in}
 \;\; L^2(0,t^*(M);\R).
 \ee
 Clearly, by the strong bang-bang property
  of $(TP)^M$
  (see Lemma~\ref{Lemma 1}), it follows that
$$|u^i_*(t)|=k_iM,\;\;\text{for a.e.}\;\;
 t\in (0,t^*(M))\;\;\text{and for all}\;\;i=1,\cdots,d.$$
 Hence, for each $i$ with $1\leq i\leq d$, it holds that
 $$\|u^i_n\|_{L^2(0,t^*(M);\R)}
 \longrightarrow\|u_*^i\|_{L^2(0,t^*(M);\R)}
 \;\;\text{as}\;\;n\rightarrow +\infty.$$
 This, along with (\ref{4.6}), leads to (\ref{4.2}) and completes
 the proof of Theorem~\ref{Theorem A}.

\bigskip
\noindent\textbf{Acknowledgements:}\;\;The author would like to appreciate Professor Gengsheng Wang for his valuable help for this paper.

\end{document}